# Improvement on Brook's theorem for $2K_2$–free Graphs

Medha Dhurandhar

**Abstract:** Here we prove that for a $2K_2$-free graph G with $\Delta(G) \geq 5$, $\chi(G) \leq \max\{\Delta-1, \omega\}$. This implies that **Borodin & Kostochka Conjecture** is true for $2K_2$-free graphs.

**Introduction:**

In 1941, Brooks' theorem states that for any connected undirected graph *G* with maximum degree $\Delta$, the ch romatic number of *G* is at most $\Delta$ unless *G* is a complete graph or an odd cycle, in which case the chromatic number is $\Delta + 1$ [1]. In 1977, **Borodin & Kostochka** conjectured that if $\Delta(G) \geq 9$, then $\chi(G) \leq \max\{\omega, \Delta-1\}$ [2]. In 1999, Reed proved the conjecture for $\Delta \geq 10^{14}$ [3]. Also D. W. Cranston and L. Rabern [4] proved it for claw-free graphs. Here we prove that $2K_2$-free graphs are $\Delta$-1-colorable. From this **Borodin & Kostochka** conjecture follows as a corollary for $2K_2$-free graphs.

**Notation:** For a graph G, V(G), E(G), $\Delta$, $\omega$, $\chi$ denote the vertex set, edge set, maximum degree, size of a maximum clique, chromatic number of G resply. For $u \in V(G)$, $N(u) = \{v \in V(G) / uv \in E(G)\}$, and $\overline{N(u)} = N(u) \cup (u)$. If $S \subseteq V$, then $<S>$ denotes the subgraph of G induced by S. If C is some coloring of G and if a vertex u of G is colored m in C, then u is called a m-vertex. Also if P is a path in G s.t. vertices on P are alternately colored say i and j, then P is called an i-j path. All graphs considered henceforth are simple. We consider here simple and undirected graphs. For terms which are not defined herein we refer to Bondy and Murty [5].

Note that $C_5$ is the only odd cycle in a $2K_2$-free graph.

**Main Result:** If $\Delta \geq 5$, and G is $2K_2$-free, then $\chi \leq \max\{\omega, \Delta-1\}$.

Proof: Let G be a smallest $2K_2$-free with $\Delta \geq 5$ and $\chi > \max\{\omega, \Delta-1\}$. Then clearly as $G \neq C_{2n+1}$ or $K_{|V(G)|}$, $\chi = \Delta > \omega$. Let $u \in V(G)$. Then $G-u \neq K_{|V(G)|-1}$ (else $\chi = \omega$). If $\Delta(G-u) \geq 5$, then by minimality $\chi(G-u) \leq \max\{\omega(G-u), \Delta(G-u)-1\}$. Clearly if $\omega(G-u) \leq \Delta(G-u)-1$, then $\chi(G-u) = \Delta(G-u)-1 \leq \Delta-1$ and otherwise $\chi(G-u) = \omega(G-u) \leq \omega < \Delta$. In any case, $\chi(G-u) \leq \Delta-1$. Also if $\Delta(G-u) < 5$, then as $G-u \neq C_{2n+1}$ (else as G is $2K_2$-free, G-u ~ $C_5$, G ~ $W_6$ and $\chi = \Delta-1$), by Brook's Theorem $\chi(G-u) \leq \Delta(G-u) < 5 \leq \Delta$. Thus always $\chi(G-u) \leq \Delta-1$ and in fact, $\chi(G-u) = \Delta-1$ and deg $v \geq \Delta-1$ $\forall$ $v \in V(G)$.

Let $Q \subseteq V(G)$ be s.t. $<Q>$ is a maximum clique in G. Let $u \in Q$ be s.t. deg $u = \max\limits_{v \in Q}$ deg v. As $\Delta > \omega$, $|N(u)-Q| \geq 1$. Let $S = \{1,..., \Delta\}$ be a $\Delta$-coloring of G s.t. u is colored $\Delta$ and vertices in Q are colored 1,.., $\omega$-1. Every vertex v of N(u) with a unique color say i has at most two vertices of the same color (else $\overline{N(v)}$ has a color say r missing. Then color v by r and u by i).     **I**

**Case 1:** deg u = $\Delta$-1 $\Rightarrow$ deg v = $\Delta$-1 $\forall$ v $\in$ Q.
Clearly every j-vertex v of Q has a unique i-vertex in N(v) where i $\neq$ j (else N(v) has some color r missing in $\overline{N(v)}$. Color v by r, u by j). W.l.g. let B $\in$ N(u)-Q be colored $\omega$ and A $\in$ Q be s.t. AB $\notin$ E(G) and A is colored 1. Let R be a component containing A s.t. vertices in R are alternately colored 1 and $\omega$. Then B $\in$ R (else alter colors in R and color u by 1). Let P = {A, D, E, B} be a 1-$\omega$ path in R. Also as $|N(u)| \geq 4$, let F, G $\in$ N(u) be colored 2, 3 resply.

Now FD or FE ∈ E(G) (else <uF∪DE> = 2K₂). W.l.g. let FD ∈ E(G). As D is the only ω-vertex of F, FB ∉ E(G). Let H be the unique 2-vertex of B ⇒ DH ∈ E(G) (else <FD∪BH> = 2K₂). Now GD ∈ E(G) (else if X ≠ G is a 3-vertex of D, then <uG∪DX> = 2K₂ and if D has no 3-vertex, then color D by 3, F by ω and u by 2) ⇒ D is the only ω-vertex of G and D has two 1-vertices and 2-vertices. If G is the only 3-vertex of D, then color D by 3, G by ω and u by 3 and if D has another 3-vertex, then $\overline{N(D)}$ has a color say r missing, color D by r, A by ω and u by 1, a contradiction in both the cases.

**Case 2:** deg u = Δ ≥ 5.
Let A, B, C, D, E ∈ N(u) and A, B, C, D be colored 1, 2, 3, 4 resply.

**Case 2.1:** ω = 2
As in **Case 1** let P = {A, F, G, B} be a 1-2 path. Then as ω = 2, C is adjacent to exactly one of F or G (else if CF, CG ∉ E(G), then <uC∪FG> = 2K₂). W.l.g. let CG ∈ E(G) ⇒ CF ∉ E(G). Now C has a 2-vertex say H (else color C by 2, u by 3) ⇒ GH ∉ E(G). But then if AH ∉ E(G), <AF∪CH> = 2K₂ and if AH ∈ E(G), <AH∪BG> = 2K₂, a contradiction in both the cases.

**Case 2.2:** ω > 2

**Case 2.2.1:** N(u)-Q has a repeat color
W.l.g. let A, B ∈ Q. Let E be colored 4 ⇒ D, E ∈ N(u)-Q.

**Case 2.2.1.1:** ∃ a vertex W ∈ Q s.t. WD, WE ∉ E(G). W.l.g. let W = A.
As before ∃ a 1-4 path P from A containing either D or E (else alter colors along P, color u by 1). W.l.g. let P = {A, F, G, D}. Then EG ∈ E(G) (else <uE∪FG> = 2K₂). Also BG ∈ E(G) (else as <uB∪FG> ≠ 2K₂, BF∈ E(G); as <AB∪DG> ≠ 2K₂, BD ∈ E(G); and as <AB∪EG> ≠ 2K₂, BE∈ E(G) contrary to **I**).

**Case 2.2.1.1.1:** BF ∉ E(G).
Then F has no 2-vertex (else if X is a 2-vertex of F, then <uB∪FX> = 2K₂) ⇒ A has one more 4-vertex say H (else color F by 2, A by 4, u by 1) and hence by I, B is its only one 2-vertex. Then GH ∈ E(G) (else <AH∪DG> = 2K₂) ⇒ ∃ color say r ∈ S s.t. G has no r-vertex as deg G ≤ Δ. Also as B has two 1-vertices it has a unique 4-vertex. W.l.g. let BE ∉ E(G). Then E has no 2-vertex (else if X is a 2-vertex of E, then <AB∪EX> = 2K₂) ⇒ BD ∈ E(G) (else as before D has no 2-vertex. Color D, E by 2 and u by 4). Also G is the only 1-vertex of E (else if X ≠ G is a 1-vertex of E, then <AB∪EX> = 2K₂) ⇒ D has two 1-vertices (else color G by r, D & E by 1, u by 4) and hence B is the only 2-vertex of D. Color D by 2, B by 4, A by 2 and u by 1, a contradiction.

**Case 2.2.1.1.2:** BF ∈ E(G).
As B has two 1-vertices, F is its only 4-vertex and BD, BE ∉ E(G) ⇒ F has another 2-vertex (else color B by 4, F by 2, u by 2). Also FC ∈ E(G) (else if F has a 3-vertex say X ≠ C, then <uC∪FX> = 2K₂ and if F has no 3-vertex, then color F by 3, B by 4 and u by 2). Further as F has two 1-vertices and 2-vertices C is its only 3-vertex (else some color r is missing in $\overline{N(F)}$. Then color F by r, B by 4, u by 2). Also G is only 1-vertex of D and E (else if X ≠ G is a 1-vertex of say D, then if FX ∉ E(G), <AF∪DX> = 2K₂ and if FX ∈ E(G), then some color r being missing in $\overline{N(F)}$, color F by r, B by 4, u by 2) ⇒ as before GC ∈ E(G) and C is the only 3-vertex of G. If C has one 4-vertex F, then color C by 4, F by 3, u by 3 and if C has one 1-vertex G, then color C by 1, G by 3, u by 3, a contradiction.

**Case 2.2.1.2:** ∀ V ∈ Q either VD or VE ∈ E(G).
Now D is non-adjacent to some vertex of Q. W.l.g. let DA ∉ E(G). Then AE ∈ E(G). Again w.l.g. let EB ∉ E(G). Then BD ∈ E(G).

**Case 2.2.1.2.1:** E is the only 4-vertex of A.

Then E has a 2-vertex say F (else color E by 2, A by 4, u by 1) and DF ∈ E(G) (else <EF∪BD> = $2K_2$).

A) If G is another 2-vertex of E, then DG ∈ E(G) (else <EG∪BD> = $2K_2$) and A is the only 1-vertex of E ⇒ B has another 1-vertex (else color A by 4, E by 1, B by 1, u by 2) and D is the only 4-vertex of B. But then as D has three 2-vertices some color r is missing in $\overline{N(D)}$. Color D by r, B by 4, u by 2, a contradiction.

B) F is the only 2-vertex of E. Then FA ∈ E(G) (else as before color 1 is missing in $\overline{N(F)}$. color F by 1, E by 2, A by 4, u by 1) and F has another 1-vertex (else color F by 1, E by 2, A by 4, u by 1). Also CF ∈ E(G) (else if F has no 3-vertex, then color F by 3, E by 2, A by 4, u by 1 and If X is a 3-vertex of F, then <uC∪FX> = $2K_2$) and C is the only 3-vertex of F (else as F has two 1-vertices and 4-vertices, some color r is missing in $\overline{N(F)}$. Color F by r, E by 2, A by 4, u by 1). Now C has two 2-vertices (else color C by 2, F by 3, u by 3). Hence C has one 1-vertex. Again CA ∉ E(G) (else color A by 4, E by 2, F by 3, C by 1, u by 3). Then as ∃ 1-3 path from A to C, B has either two 1-vertices or 3-vertices and hence D is the only 4-vertex of B and as before D has a 1-vertex say G, GC ∈ E(G) and C (G) is the only 3-vertex (1-vertex) of G (C). Color C by 1, G by 3, u by 3, a contradiction.

**Case 2.2.1.2.2:** Every W ∈ Q has two 4-vertices.

Then C ∈ Q (else if CW ∉ E(G) for some W ∈ Q, where W has color r, then as before A has either two r-vertices or 3-vertices) ⇒ C has two 4-vertices and CD or CE ∈ E(G) (else we get **Case 2.2.1.1**). W.l.g. let CE ∈ E(G). Also as A, B, C have two 4-vertices each, clearly B, C (C, A; A, B) are the only 2, 3 (3, 1; 1, 2) vertices of A (B; C).

A) CD ∉ E(G).

Then as seen earlier D has no 3-vertex and E has two 1-vertices (else color D by 3, A by 2, B by 1, E by 1, u by 4). Similarly E has two 3-vertices and some color r is missing in $\overline{N(E)}$. Color D by 3, E by r, u by 4, a contradiction.

B) CD ∈ E(G).

    **B.1**    If A is the only 1-vertex of E, then clearly D has no 1-vertex. Hence color D by 1, C by 4, E by 1, A by 3, u by 4 ⇒ color 1 is repeated in N(u)-Q, A has only one 1-vertex E and we get **Case 2.2.1.2.1**, a contradiction.

    **B.2**    Hence E (D) has two 1 (2)-vertices viz. A, F (B, G). Also DF, EG ∈ E(G) and FG ∉ E(G) (else <AB∪FG> = $2K_2$). As CF, CG ∉ E(G) F, G have no 3-vertex. Hence color F, G by 3. If A (B) is the only 1-vertex (2-vertex) of E (D), then we get **B.1**. Hence E (D) has one more 1-vertex (2-vertex) and some color r (s) is missing in $\overline{N(E)}$, ($\overline{N(D)}$). Color D by s, E by r, u by 4, a contradiction.

**Case 2.2.2:** A ∈ Q and E ∈ N(u)-Q have a common color 1.

**Case 2.2.2.1:** ω = 3.

Then <B, C, D> = $\overline{K_3}$ (else this case reduces to **Case 2.2.1**). Now B has either one 3-vertex or 4-vertex (else some color r is missing in $\overline{N(B)}$. Color B by r, u by 2). W.l.g. let B have one 3-vertex. As before let P = {B, F, G, C} be a 2-3 path. Then FD ∈ E(G) (else if F has a 4-vertex X ≠ D, then <uD∪FX>= $2K_2$ and if F has no 4-vertex, then color F by 4, B by 3, u by 2). Now C has a 4-vertex say H (else color C by 4, u by 3). Then FH ∈ E(G) (else <CH∪DF> = $2K_2$). Now F has two 2-vertices and 4-vertices, hence G is the only 2-vertex of C (else if X ≠ G is another 2-vertex of C, then either FX ∉ E(G) and <BF∪CX> = $2K_2$ or FX ∈ E(G) and ∃ color r missing in $\overline{N(F)}$ in which case color F by r, B by 3, u by 2). Similarly H is the only 4-vertex of C. Also DG ∈ E(G) (else if X ≠ G is a 2-vertex of E, then <CG∪DX> = $2K_2$). Similarly BH ∈ E(G) ⇒ GH ∈ E(G) (else <BH∪DG> =

$2K_2$). As $\omega = 3$, either AF or AH $\notin$ E(G). W.l.g. let AF $\notin$ E(G). Then AG (AH) $\in$ E(G) (else <uA∪FG> (<uA∪HG>) = $2K_2$). Also A is the only 1-vertex of G (H) (else as before G (H) has a color say r missing, then color G (H) by r, C by 2 (4), u by 3) $\Rightarrow$ <uE∪GH> = $2K_2$, a contradiction.

**Case 2.2.2.2:** $\omega > 3$. Let $C \in Q$.
W.l.g. let DB $\notin$ E(G) (else if DW $\in$ E(G) $\forall$ W $\in$ Q-A, then with Q' = Q-A+D we get Case 2.2.1). As before if P = {B, F, G, E} is a 2-4 path, then A, C have either two 2-vertices or 4-verices. Thus A is the only 1-vertex of C (else C has a color say r missing in which case color C by r, u by 3) $\Rightarrow$ C is the only 3- vertex of A (else A has a color say s missing in which case color A by s, C by 1, u by 3).   **II**
Then CE $\notin$ E(G) and E has no 3-vertex (else if X is a 3-vertex of E, then <AC∪EX>= $2K_2$).

**Case 2.2.2.2.1:** $\exists$ W $\in$ Q-{A, B} s.t. DW $\notin$ E(G). W.l.g. let W = C.
Then as before A (B) is the only 1-vertex (2-vertex) of B (A) and AG $\notin$ E(G) $\Rightarrow$ AF, AD $\in$ E(G) (else <uA∪FG>= $2K_2$ or <AB∪EG>= $2K_2$). Also G has no 1-vertex (else if X is a 1-vertex of G, then <AB∪GX>= $2K_2$) $\Rightarrow$ D has another 2-vertex X $\neq$ G (else color G by 1, D by 2, u by 4). Hence A is the only 1-vertex of D. Color A by 2, B by 1, D by 1, u by 4, a contradiction.

**Case 2.2.2.2.2:** DC $\in$ E(G) $\forall$ C $\in$ Q-{A, B}.

A. AD $\in$ E(G)
   Then G is symmetric about B and D. As before let P = {B, F, G, D} be a 2-4 path and w.l.g. let CF $\in$ E(G). As C has two 4-vertices, C has only one 2-vertex B and hence CG $\notin$ E(G). Thus G has no 3-vertex (else if X is a 3-vertex of G, then <uC∪GX>= $2K_2$) $\Rightarrow$ D has another 2-vertex say H $\neq$ G (else color G by 3, D by 2, u by 4) and FH $\in$ E(G) (else <BF∪DH>= $2K_2$). Then C is the only 3-vertex of D (else some color r is missing in $\overline{N(D)}$. Color D by r, u by 4). Now B has another 3-vertex (else color B by 3, C by 2, D by 3, u by 4) $\Rightarrow$ F is the only 4-vertex of B and as F has three 2-vertices C is the only 3-vertex of F (else some color r is missing in $\overline{N(F)}$. Color F by r, B by 4, u by 2). But then color C by 2, B by 4, F and u by 3, a contradiction.

B. AD $\notin$ E(G)
   Then by **II**, D has another 3-vertex say H $\neq$ C (else color A by 3, C by 1, D by 3, u by 4) and hence G is the only 2-vertex of D. Now BH $\in$ E(G) (else <AB∪DH>= $2K_2$) $\Rightarrow$ F is the only 4-vertex of B. Again CF $\notin$ E(G) or CG $\notin$ E(G) (else some color r is missing in $\overline{N(C)}$. Color C by r, u by 3). Now if CF (CG) $\notin$ E(G), then F (G) has no 3-vertex, hence color F (G) by 3, B by 4 (D by 2), u by 2 (4) , a contradiction.

   This proves the theorem.

**Corollary:** Borodin and Kostochka conjecture is true for $2K_2$-free graphs.

**Result 2:** Let G be $2K_2$-free. Then either $\chi \leq \max\{\omega, \Delta-1\}$ or
1) $\Delta \leq 3$, $\omega = 2$, $\chi = 3$ and G has an odd cycle induced or
2) $\Delta = 4 = \chi$, $\omega = 3$ and G is one of the eight graphs shown in **Figure 1** below.

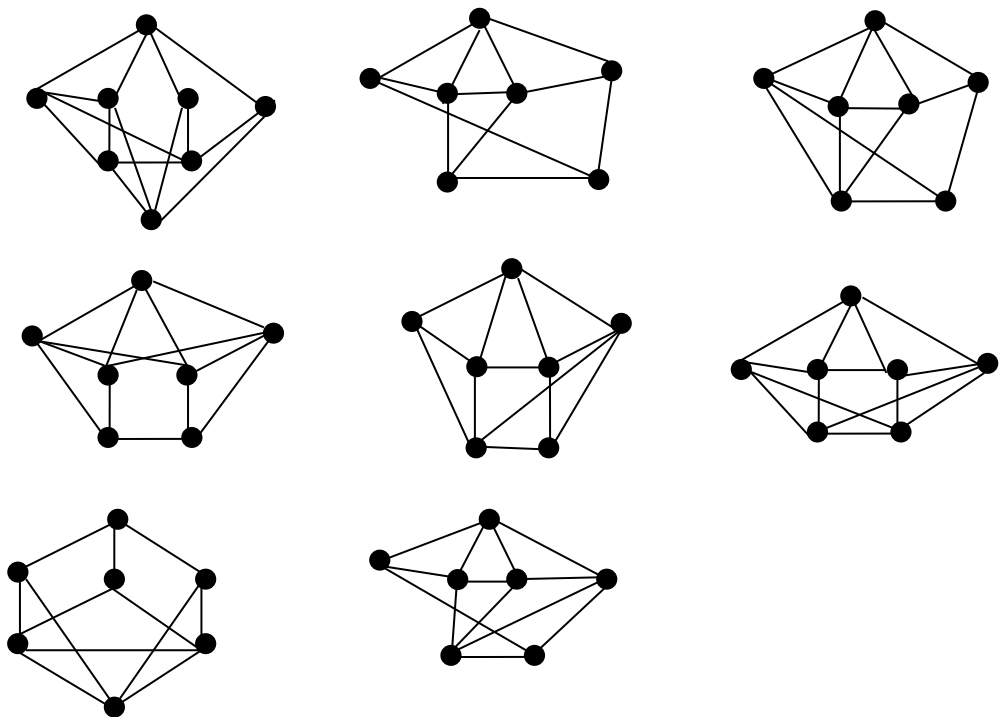

**Figure 1**

**References**

[1] Brooks, R. L., "On colouring the nodes of a network", Proc. Cambridge Philosophical Society, Math. Phys. Sci., 37 (1941), 194–197

[2] O. V. Borodin and A. V. Kostochka, "On an upper bound of a graph's chromatic number, depending on the graph's degree and density", JCTB 23 (1977), 247--250.

[3] B. A. Reed, "A strengthening of Brooks' Theorem", J. Comb. Theory Ser. B, 76 (1999), 136–149.

[4] D. W. Cranston and L. Rabern, "Coloring claw-free graphs with ∆-1 colors" *SIAM J. Discrete Math.*, 27(1) (1999), 534–549.

[5] J.A. Bondy and U.S.R. Murty. Graph Theory, volume 244 of Graduate Text in Mathematics. Springer, 2008.